\documentclass[number,citesort,seceqn,dvips]{arxbj}
\usepackage{cursive} %uplgreek}
\usepackage{graphicx}

\aid{0}
\volume{16}
\issue{4}
\pubyear{2010}
\firstpage{1208}
\lastpage{1223}
\doi{10.3150/10-BEJ248}

\makeatletter
\newremark{example}{Example}
\renewcommand{\citep}[1]{\cite{#1}}

\makeatother

\begin{document}
\begin{frontmatter}

\title{Second order ancillary: A differential view from continuity}
\runtitle{Second order ancillary}

\begin{aug}
\author[a]{\fnms{Ailana M.} \snm{Fraser}\corref{}\thanksref{a}\ead[label=e1]{afraser@math.ubc.ca}},
\author[b]{\fnms{D.A.S.} \snm{Fraser}\thanksref{b}\ead[label=e2]{dfraser@utstat.toronto.edu}} \and
\author[c]{\fnms{Ana-Maria}~\snm{Staicu}\thanksref{c}\ead[label=e3]{staicu@stat.ncsu.edu}}
\runauthor{A.M. Fraser, D.A.S. Fraser and A.-M. Staicu}
\address[a]{Department of Mathematics, University of British Columbia,
Vancouver, Canada V6T 1Z2.\\ \printead{e1}}
\address[b]{Department of Statistics, University of Toronto,
Toronto, Canada M5S 3G3.\\ \printead{e2}}
\address[c]{Department of Statistics, North Carolina State University,
Raleigh, NC 27695, USA.\\ \printead{e3}}
\end{aug}

% HISTORY:
\received{\smonth{1} \syear{2009}}
\revised{\smonth{12} \syear{2009}}

% ABSTRACT
%
\begin{abstract}
Second order approximate ancillaries have evolved as the primary
ingredient for recent likelihood development in statistical inference.
This uses quantile functions rather than the equivalent distribution
functions, and the intrinsic ancillary contour is given explicitly as
the plug-in estimate of the vector quantile function.
The derivation uses a Taylor expansion of the full quantile function,
and the linear term gives a tangent to the observed ancillary contour.
For the scalar parameter case, there is a vector
field that integrates to give the ancillary contours, but for the
vector case,
there are multiple vector fields and the Frobenius conditions for
mutual consistency
may not hold. We demonstrate, however, that the conditions hold in a
restricted way and that this verifies the
second order ancillary contours in moderate deviations.
The methodology
can generate an appropriate exact ancillary when such exists or an
approximate ancillary for the numerical or Monte Carlo calculation of
$p$-values and confidence quantiles. Examples are given, including
nonlinear regression and several enigmatic examples from the literature.
\end{abstract}

% KEYWORDS
%

\end{frontmatter}

%s1 ###
\section{Introduction}\label{intro}

Ancillaries are loved or hated, accepted or rejected, but typically
ignored. Recent approximate ancillary methods (e.g., \citep{reid}) give a decomposition of the sample space rather than providing
statistics on the sample space (e.g., \citep{cox1,mcc1}). As a result, continuity gives the contour along which the
variable directly measures the parameter and then gives the subcontour
that provides measurement of a parameter of interest. This, in turn,
enables the high accuracy of cumulant generating function
approximations \citep{dan,bar1} to extend to cover a wide
generality of statistical models.

Ancillaries initially arose (see \citep{fis1}) to examine the accuracy
of the maximum likelihood estimate, then (see \citep{fis2}) to
calibrate the loss of information in the use of the maximum likelihood
estimate and then (see \citep{fis3}) to develop a key instance
involving the configuration statistic. The configuration of a sample
arises naturally in the context of sampling a location-scale model,
where a standardized coordinate $z = (y - \mu)/\sigma$ has a fixed and
known error distribution $g(z)$: the $i$th coordinate of the response
thus has $f(y_i; \mu, \sigma) = \sigma^{-1} g\{(y_i - \mu)/\sigma\}$.
The configuration $a(y)$ of the sample is the plug-in estimate of the
standardized residual,
%
%e1.1 ###
\begin{equation}
\label{hat}
a(y) = \hat z =  \biggl({y_1 - \hat\mu\over\hat\sigma}, \ldots, {y_n -
\hat\mu\over\hat\sigma}  \biggr)^\prime,
\end{equation}
where $(\hat\mu, \hat\sigma)$ is the maximum likelihood value for $(\mu
, \sigma)$ or is some location-scale equivalent. Clearly, the
distribution of $\hat z$ is free of $\mu$ and $\sigma$ as the
substitution $y_i = \mu+ \sigma z_i$ in (\ref{hat}) leads to the
cancellation of dependence on $\mu$ and $\sigma$. This supports a
common definition for an ancillary statistic $a(y)$, that it has a
parameter-free distribution; other conditions are often added to seek
sensible results.

More generally, the observed value of an ancillary identifies a sample
space contour along which parameter change modifies the model, thus
yielding the conditional model on the observed contour as the
appropriate model for the data. The ancillary method is to use directly
this conditional model identified by the data.

One approach to statistical inference is to use only the observed
likelihood function $L^0(\theta) = L(\theta; y^0)$ from the model $f(y;
\theta)$ with observed data $y^0$. Inference can then be based on some
simple characteristic of that likelihood. Alternatively, a weight
function $w(\theta)$ can be applied and the composite $w(\theta)L(\theta
)$ treated as a distribution describing the unknown $\theta$; this
leads to a rich methodology for exploring data, usually, but
unfortunately, promoted solely within the Bayesian framework.

A more incisive approach derives from an enriched model which is often
available and appropriate. While the commonly cited model is just a set
of probability distributions on the sample space, an enriched model can
specifically include continuity of the model density function and
continuity of coordinate distribution functions. An approach that
builds on these enrichments can then, for example, examine the observed
data $y^0$ in relation to other data points that have a similar shape
of likelihood and are thus comparable, and can do even more. For the
location-scale model, such points are identified by the configuration
statistic; then, accordingly, the model for inference would be $f\{y
\mid
a(y) = a^0 ;\theta\}$, where $a(y)$ is the configuration ancillary.

Exact ancillaries as just described are rather rare and seem limited to
location-type models and simple variants. However, extensions that use
approximate ancillaries
(e.g., \citep{fra5,fra9}) have recently been broadly fruitful,
providing approximation in an asymptotic sense. Technical issues can
arise with approximate values for an increasing number of coordinates,
but these can be managed by using ancillary contours rather than
statistics; thus, for a circle, we use explicitly a contour
$A = \{(x,y)= (a^{1/2}\cos t, a^{1/2}\sin t) \dvt  t \mbox{ in }[0, 2\curpi)\}$
rather than using implicitly a statistic $x^2 + y^2 = a$.

We now assume independent coordinate distribution functions that are
continuously differentiable with respect to the variable and the
parameter; extensions will be discussed separately. Then, rather than
working directly with a coordinate distribution function $u_i = F_i
(y_i; \theta)$, we will use the inverse, the quantile function $y_i =
y_i(u_i; \theta)$ which presents a data value $y_i$ in terms of a
corresponding $p$-value $u_i$. For additional advantage, we could use a
scoring variable $x$ in place of the $p$-value, for example, $x = \Phi
^{-1}(u)$ or $x = F^{-1}(u; \theta_0)$, where $\Phi(\cdot)$ is the
standard Normal distribution function. We can then write $y = y(x;
\theta)$, where a coordinate $y_i$ is presented in terms of the
corresponding scoring variable $x_i$.

For the full response variable, let
$
y = y(x; \theta) =
\{y_1(x_1; \theta),\ldots, y_n(x_n; \theta) \}'
$
be the quantile vector expressing $y$ in terms of the reference or
scoring variable $x$ with its given distribution: the quantile vector
records how parameter change affects the response variable and its
distribution, as prescribed by the continuity of the coordinate
distribution functions.

For an observed data point $y^0$, a convenient reference value $\hat
x^0$ or the fitted $p$-value vector is obtained by solving the equation
$ y^0 = y(x; \hat\theta^0)$
for $x$, where $\hat\theta^0$ is the observed maximum likelihood value;
for this, we assume regularity and asymptotic properties for the
statistical model.
The contour of the second order ancillary through the observed data
point as developed in this paper is then given
as the trajectory of the reference value,
%
%e1.2 ###
\begin{equation}\label{obs}
A^0= \{ y(\hat x^0;t) \dvt  t \mbox{ in } \mathbb{R}^p\},
\end{equation}
to second order under parameter change,
where $p$ here is the dimension of the parameter. A~sample space point
on this contour has, to second order, the same estimated $p$-value
vector as the observed data point and special
properties for the contours are available to second order.

The choice of the reference variable with given data has no effect on
the contour: the reference variable could be Uniform, as with the
$p$-value; or, it could be the response distribution itself for some
choice of the parameter, say $\theta_0$.

For the location-scale example mentioned earlier, we have the
coordinate quantile function $y_i = \mu+ \sigma z_i$, where $z_i$ has
the distribution $g(z)$. The vector quantile function is
%
%e1.3 ###
\begin{equation} \label{quan}
y(z; \mu, \sigma) = \mu1+\sigma z ,
\end{equation}
where $1 = (1, \ldots, 1)^\prime$ is the `one vector.' With the data
point $y^0$, we then have the fitted $\hat z^0 =(y^0 - \hat\mu^0 1)/
\hat\sigma^0$. The observed ancillary contour to second order is then
obtained from (\ref{obs}) by substituting $\hat z^0$ in the quantile
(\ref{quan}):
%
%e1.4 ###
\begin{equation} \label{lsanc}
A^0 = \{y(\hat z^0; t)\} = \{m1 + s\hat z^0; (m, s) \mbox{ in } \mathbb{R}\times\mathbb{R}^+\} =\mathcal{L}^+(1; \hat z^0)
\end{equation}
with positive coefficient for the second vector.
This is the familiar exact ancillary contour $a(y)=a^0$ from (\ref{hat}).

An advantage of the vector quantile function in the context of the
enriched model mentioned above
is that it allows us to examine how parameter change modifies the
distribution and thus how it moves data points as a direct expression
of the explicit continuity.
In this sense, we define the velocity vector or vectors as
$v(x; \theta) = (\partial/ \partial\theta)y(x; \theta) =\partial y/
\partial\theta$.
In the scalar $\theta$ case, this is a vector recording the direction
of movement of a point $y$ under $\theta$ change; in the vector $\theta
$ case, it is a $1 \times p$ array of such vectors in $\mathbb{R}^n$,
$V(x; \theta) = \{v_1(x_1; \theta), \ldots, v_p(x_p; \theta)\},
$
recording the separate effects from the parameter coordinates $\theta
_1, \ldots, \theta_p$.
For the location-scale example, the velocity array is
$V(z; \mu, \sigma) = (1, z),$
which can be viewed as a $1 \times2$ array of vectors in $\mathbb{R}^n$.

The ancillary contour can then be presented using a Taylor series about
$y^0$ with coefficients given by the velocity and acceleration $V$ and $W$.
For the location-scale example,
the related acceleration vectors are equal to zero.

For more insight, consider the general scalar $\theta$ case and the
velocity vector $v(x; \hat\theta^0)$. For a typical coordinate, this
gives the change $\mathrm{d}y=v(x; \hat\theta^0)\,\mathrm{d}\theta$ in the variable as
produced by a small change $\mathrm{d}\theta$ at $\hat\theta^0$. A re-expression
of the coordinate variable can make these increments equal and produce
a location model; the product of these location models is a full
location model $g(y_1 - \theta, \ldots, y_n - \theta)$ that precisely
agrees with the initial model to first derivative at $\theta= \hat
\theta^0$ (see \citep{fra7,and}). This location model then, in
turn, determines a full location ancillary with configuration $a(y) =
(y_1 - \bar y, \ldots, y_n - \bar y)$. For the original model, this
configuration statistic has first-derivative ancillarity at $\theta=
\hat\theta^0$ and is thus a first order approximate ancillary; the
tangent to the contour at the data point is just the vector $v(\hat
x^0; \hat\theta^0)$. Also this contour can be modified to give second
order ancillarity.

%The quantile function allows us to examine directly how the parameter
%moves the distribution; the vector
%$$
%v(y; \theta) = {dy \over d\theta} = {\partial y(z; \theta) \over
%$$
%records how parameter change at $\theta$ affects the distribution as
%identified by its $p$-value or scoring variable. If we then examine
%this for $\theta= \hat\theta^0$ we are led in the scaler $\theta$
%case to a location model $g\{y_1 - \theta, \ldots, y_n -\theta\}$ that
%precisely agrees with the given model for $\theta$ near $\hat\theta^0$
%and then has correspondingly the full location ancillary available
%(Fraser \& Reid, 2002; Andrews, Fraser \& Wong, 2005); this is a
%primary input for the ancillary developed here.

In a somewhat different way, the velocity vector $v(y^0; \theta)$ at
the data point $y^0$ gives information as to how data change at $y^0$
relates to parameter change at various $\theta$ values of interest.
This allows us to examine how a sample space direction at the data
point relates to estimated $p$-value and local likelihood function
shape at various $\theta$
values; this, in turn, leads to quite general default priors for
Bayesian analysis (see \citep{fra8}).

In the presence of a cumulant generating function, the saddle-point
method has produced highly accurate third order approximations for
density functions (see \citep{dan}) and for distribution functions (see
\citep{lug}). Such approximations are available in the presence of
exact ancillaries \citep{bar1} and extend widely in the presence of
approximate ancillaries (see \citep{fra5}). For third order accuracy,
only second order approximate ancillaries are needed, and for such
ancillaries, only the tangents to the ancillary contour at the data
point are needed (see \citep{fra5,fra6}). With this as our
imperative, we develop the second order ancillary for statistical inference.

Tangent vectors to an ancillary at a data point
give information as mentioned above concerning a location model approximation
at the data point. For a scalar parameter, these provide a vector field
and integrate quite generally to give a unique approximate ancillary to
second order
accuracy. The resulting
conditional model then provides definitive $p$-values by available
theory; see, for example,
\citep{fra9}.
For a vector parameter, however, the multiple vector fields may not satisfy
the Frobenius conditions for integrability and thus may not define a function.

Under mild conditions, however, we show that such tangent vectors do generate
a surface to second order without the Frobenius conditions holding.
We show this in several steps. First, we obtain the coordinate quantile
functions
$
y_i = y_i(x_i; \theta).
$
Second, we Taylor series expand the full vector quantile
$y=(y_1,\ldots,y_n)$ in terms
of the full reference variable $x=(x_1,\ldots,x_n)$ and the parameter
$\theta=(\theta_1,\ldots,\theta_p)$ about data-based values, appropriately
re-expressing coordinates and working to second order.
Third, we
show that this generates a partition with second order ancillary properties
and the usual tangent vectors. The seeming need for the full Frobenius
conditions is bypassed
by finding that two integration routes need not converge to each other, but
do remain on the same contour, calculating, of course, to second order.

%Let $y = y(x;\theta)$ be the full quantile vector; then with observed
%data $y^0$ let
%$\hat x^0$ be the
%fitted reference value which satisfies
%$$
%y^0=y(\hat x^0; \hat\theta^0)
%$$
%where $\hat\theta^0$ is the observed maximum likelihood value.
%The contour of the second order ancillary through the observed data
%$y^0$
%is then given as
%A^0= \{ y(\hat x^0;t) : t {\in} R^p\}.
%A sample space point on this contour has for some $\theta$ value the
%same vector of
%coordinate $p$-values as the observed data $y^0$ has at its maximum
%likelihood
%parameter value; the contour thus represents constant $p$-value
%information.

%For an exact ancillary that has full continuity in the model,
%this reproduces the local form of that
%ancillary contour through the data point. And more generally it gives
%a local approximate ancillary that conforms to continuity in the
%model; and when
%expanded in
%a Taylor expansion produces what can be called velocity and
%acceleration vectors for
%the ancillary contour; these vectors describe the first and second
%order characteristics
%of the surface.

This construction of an approximate ancillary is illustrated in Section~\ref{sec2} using the familiar example, the Normal-on-the
circle from \citep{fis4}; see also \citep{cox2,bar2,fra7,fra3}. The example, of course, does have an exact
ancillary and the present procedure gives an approximation to that ancillary.
In Section~\ref{examples}, we consider various examples that have exact and
approximate ancillaries, and then in Sections~\ref{sec4} and~\ref{sec5}, we present the
supporting theory.
In particular, in Section~\ref{sec4}, we develop notation for a $p$-dimensional
contour in $\mathbb{R}^n$,
$
A= \{y(x_0;t)\dvt t \mbox{ in } \mathbb{R}^p\},
$
and use velocity and acceleration vectors to present a Taylor series
with respect to $t$. Then, in Section~\ref{sec5}, we consider a regular
statistical model with asymptotic properties and
use the notation from Section~\ref{sec4} to develop the second order ancillary contour
through an observed data point $y^0$. The re-expression
of individual coordinates, both of the variable and the parameter,
plays an essential role in the development; an asymptotic analysis is
used to establish the second order approximate ancillarity.
Section~\ref{sec6} contains some discussion.

\section{Normal-on-the-circle}\label{sec2}

We illustrate the second order approximate ancillary with a simple
nonlinear regression model,
the \textit{Normal-on-the-circle} example (see \citep{fis4}). The model
has a well-known exact ancillary.
Let $y = (y_1, y_2)^\prime$ be Normal on the plane with
mean $(\rho\cos\theta, \rho\sin\theta)'$ and variance matrix
$I/n$ with $\rho$ known. The mean is
on a circle of fixed radius $\rho$ and the distribution has rotationally
symmetric error with variances $n^{-1}$, suggesting an
antecedent sample size $n$ for an asymptotic approach. The full
$n$-dimensional case is examined as Example~\ref{ex2} in
Section~\ref{examples} and the present case derives by routine conditioning.

The distribution is a unit probability mass centered at $(\rho\cos\theta, \rho\sin\theta)'$ on the circle with radius $\rho$. If
rotations about the origin are applied to $(y_1, y_2)'$, then the
probability mass rotates about the origin, the mean moves on the circle
with radius $\rho$ and an element of
probability at a distance $r$ from the origin moves on a circle of radius $r$.
The fact that the rotations move
probability along circles but not between circles of course implies
that probability on any circle about the origin
remains constant: probability flows on the ancillary contours.
Accordingly, we have that the radial distance $r = (y_1^2
+ y_2^2)^{1/2}$ has a fixed $\theta$-free distribution and is thus ancillary.

The statistic $r(y)$ is the Fisher exact ancillary for this problem and
Fisher recommended
that inference be based on the conditional model, given the observed
ancillary contour. This conditional approach has a long but uneven
history; \citep{fra4} provides an overview and \citep{fra10} offer
links with asymptotic theory.
We develop the approximate second order ancillary and examine how it
relates to the Fisher exact ancillary.

The model for the Normal-on-the-circle has independent coordinates, so
we can invert the coordinate
distribution functions and obtain the vector
quantile function,
\[
\pmatrix{
y_1\cr
y_2
}
=\rho
\pmatrix{
\cos \theta\cr
\sin \theta
}
+
\pmatrix{
x_1\cr
x_2
}
,
\]
where the $x_i=\Phi^{-1}(u_i)/n^{1/2}$ are independent normal variables
with means $0$ and variances $n^{-1}$, and $\Phi$ is the standard
Normal distribution function.
We now examine the second order ancillary contour $A^0$ given by (\ref{obs}).

\begin{figure}

\includegraphics{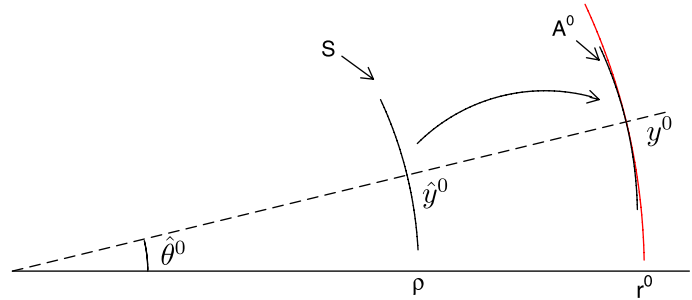}

\caption{The regression surface $S$ is a circle of radius $R$; the
local contour of
the approximate ancillary $A^0$ is a circle segment of $S$ moved from
$\hat y^0$ to $y^0$;
the exact ancillary contour is a circle segment of radius $r^0$ through
the data point $y^0$.}
\label{circle}
\end{figure}

Let $y^0 = (y_1^0, y_2^0)' = (r^0\cos a^0, r^0\sin a^0)$ be the
observed data point where $r^0$, $a^0$ are the corresponding polar
coordinates; see Figure~\ref{circle}. For this simple nonlinear normal
regression model, $\hat\theta^0 = a^0$ is the angular direction of the
data point. The fitted reference value $\hat x^0$ is the solution of
the equation
$y^0 = y(x; \hat\theta^0) =\rho(\cos a^0, \sin a^0)' + (x_1, x_2)$,
giving $\hat x^0 = (\hat x_1^0, \hat x_2^0)'= y^0 - \rho(\cos a^0,
\sin a^0)' =y^0 - \hat y^0$,
where $\hat y^0=\rho(\cos a^0, \sin a^0)'$ is the fitted value, which
is the projection of the data point $y^0$ onto the circle. The observed
ancillary contour is then
\[
A^0 =  \left\{\rho
\pmatrix{
\cos\theta\cr
\sin\theta
}
+ y^0 - \hat y^0\dvt  \theta \mbox{ near } a^0\right \}
= y^0 -\hat y^0+  \left\{\rho
\pmatrix{
\cos(a^0 + t)\cr
\sin (a^0 + t)
}
\dvt  t \mbox{ near } 0\right \}.
\]
Figure~\ref{circle} shows that $A^0 = \{y(\hat x^0; t)\dvt  t \mbox{ near }a^0\}$ is a translation, as shown by the arrow of a segment $S$ of the
solution contour, from the fitted point $\hat y^0$ to the data point
$y^0$.
%
%f1 ###

The second order ancillary segment at $y^0$ does not lie on the exact
ancillary surface $r(y_1, y_2)=r^0$. The tangent vector at the data
point $y^0$ is
$v = (\partial y / \partial t)|_{t = a^0} =
(-\rho\sin a^0, \rho\cos a^0)'$, which
is the same as the tangent vector for the exact ancillary and which
agrees with
the usual tangent vector $v$ (see \citep{fra9}). However, the
acceleration vector
is $w =(\partial^2 / \partial t^2)y|_{t = a^0} = (-\rho\sin a^0,
-\rho\cos a^0)'$,
which differs slightly from that for the exact ancillary: the
approximation has radius of curvature $\rho$, as opposed to $r^0$ for
the exact, but the difference in moderate deviations about $y^0$ can be
seen to be small and is second order.

The second order ancillary contour through $y^0$ can also be expressed
in a Taylor series as $A^0=\{y^0 + tv+wt^2/2\dvt  t \mbox{ near } 0\}$;
here, the acceleration vector $w$ is
orthogonal to the velocity vector $v$.
Similar results hold in wide generality when $y$ has dimension $n$ and
$ \theta$ has dimension $p$; further examples are discussed in the next
section and the general development follows in Sections~\ref{sec4} and~\ref{sec5}.

%s3 ###
\section{Some examples}\label{examples}

%ex1
\begin{example}[(Nonlinear regression, $\sigma_0$ known)]\label{ex1}
Consider a nonlinear regression model $y=\eta(\theta) + x$ in $\mathbb{R}^n$,
where the error $x$
is $\operatorname{Normal}(0;\sigma^2_0 I)$ and the regression or solution surface $S=\{
\eta(\theta)\}$ is smooth with
parameter $\theta$ of dimension, say, $r$. For given data point $y^0$,
let $\hat\theta^0$ be the
maximum likelihood value. The fitted value is then $\hat y^0=\eta(\hat
\theta^0)$ and the fitted reference value is $\hat x^0=y^0 - \eta(\hat
\theta^0) =y^0-\hat y^0$. The model as presented is already in quantile
form; accordingly,
$V=(\partial\eta/ \partial\theta)|_{\hat\theta^0}, W=(\partial^2\eta
/\partial\theta^2)|_{\hat\theta^0}$
are the observed velocity and acceleration arrays, respectively, and
the approximate ancillary
contour at the data point $y^0$ is
$
A^0=\{y^0 + Vt+ t'Wt/2+\cdots\dvt t \mbox{ in } \mathbb{R}^r\},
$
which is just a $y^0-\hat y^0$ translation of the solution surface
$
S=\{\hat y^0 + Vt+ t'Wt/2+\cdots\dvt t \mbox{ in }  \mathbb{R}^r\}.
$
For this, we use matrix multiplication to linearly combine the elements
in the arrays $V$ and $W$.
\end{example}

%ex2
\begin{example}[(Nonlinear regression, circle case)]\label{ex2}
As a special case, consider the regression model where the solution surface
$S=\{\eta(\theta)\}$ is a circle of radius $\rho$ about the origin;
this is the full-dimension version
of the example in Section~\ref{sec2}. For notation, let $C=(c_1,\ldots, c_n)$ be an
orthonormal
basis with vectors $c_1, c_2$ defining the plane that includes $S$. Then
$\tilde y=C'y$ provides rotated coordinates and $\tilde\eta(\theta
)=C'\eta(\theta)
=(\rho\cos\theta, \rho\sin\theta,0,\ldots,0)$ gives
the solution surface in the new coordinates.

There is an exact ancillary given by $r=(\tilde y_1^2+\tilde
y_2^2)^{1/2}$ and $(\tilde
y_3,\ldots, \tilde y_n)$; the corresponding ancillary contour through
$\tilde y^0$ is a
circle
of radius $r^0$ through the data point $y^0$ and lying in
the plane $\tilde y_3=\tilde y_3^0,\ldots,\tilde y_n=\tilde y_n^0$. The
approximate
ancillary contour is a segment of a circle of radius $\rho$ through the
data point $y^0$ and
lying in
the same plane. This directly agrees with the simple
Normal-on-the-circle example
of Section~\ref{sec2}.

For the nonlinear regression model, Severini (\citep{sev}, page 216)
proposes an approximate ancillary by using the obvious pivot $y-\eta
(\theta)$ with the plug-in maximum likelihood value
$\theta=\hat\theta$; we show that
this gives a statistic $ A(y)=y-\eta(\hat\theta)$ that can be
misleading. In the rotated coordinates, the statistic $ A(y)$ becomes
\begin{eqnarray*}
\tilde A(y)&=&(r\cos\hat\theta, r\sin\hat\theta,\tilde y_3,\ldots,\tilde
y_n)'-(\rho\cos\hat\theta, \rho\sin\hat\theta,0,\ldots,0)'\\
&=& \{(r-\rho)\cos\hat\theta, (r-\rho)\sin\hat\theta,\tilde y_3,\ldots
,\tilde
y_n\}',
\end{eqnarray*}
which has observed value $\tilde A^0 =\{(r^0-\rho)\cos\hat\theta^0,
(r^0-\rho)\sin\hat\theta^0,\tilde
y_3^0,\ldots,\tilde y_n^0\}'.$

If we now set the proposed ancillary equal to its observed value,
$\tilde A=\tilde A^0$, we obtain $\tilde y_3=\tilde y_3^0,\ldots, \tilde
y_n= \tilde
y_n^0$
and also obtain $r=r^0$ and $\hat\theta=\hat\theta^0$. Together, these
say that
$y=y^0$, and
thus that the proposed approximate ancillary is exactly equivalent to the
original response variable, which is clearly not ancillary. Severini
does note
``\ldots it does not necessarily follow that $a$ is a second-order
ancillary statistic since
the dimension of $a$ increases with $n$.'' The consequences of using the plug-in
$\hat\theta$ in the pivot are somewhat more serious: the plug-in
pivotal approach for this example does
not give an approximate ancillary.
\end{example}
%ex3
\begin{example}[(Nonlinear regression, $\sigma$ unknown)]\label{ex3}
Consider a nonlinear regression model $y=\eta(\theta) + \sigma z$ in
$\mathbb{R}^n$, where
the error $z$ is $\operatorname{Normal}(0;I)$ and the solution surface $S=\{\eta
(\theta)\}$ is
smooth with surface dimension $r$ (see \citep{fra11}). Let $y^0$ be the
observed data point and
$(\hat\theta^0, \hat\sigma^0)$ be the corresponding maximum likelihood
value. We then have the fitted regression $\hat y^0$, the fitted residual
$\hat x^0 = y^0 -\hat y^0$, and the fitted reference value $\hat
z^0=\hat x^0/\hat\sigma^0$
which is just the standardized
residual.

Simple calculation gives the velocity and acceleration arrays
\[
\bar V=(V   \hat z^0),\qquad
\bar W=
\pmatrix{
W & 0 \cr
0 & 0
}
\]
using $V$ and $W$ from Example~\ref{ex1}. The approximate ancillary contour at
the data point
$y^0$ is then
\begin{eqnarray*}
\tilde A^0 &=& \{y^0 + VT +t'Wt/2 + \cdots+s\hat z^0\dvt  t  \mbox{ in } \mathbb{R}^r, s \mbox{ in } \mathbb{R}^+ \} \\
&=& \{\eta(t) + s \hat z^0 \dvt  t \mbox{ in } \mathbb{R}^r, s \mbox{ in } \mathbb{R}^+\} \\
&=& A^0 + \mathcal{L}^+(\hat z^0) ,
\end{eqnarray*}
where $A^0$ is as in Example~\ref{ex1}. This is the solution surface from
Example~\ref{ex1}, translated from $\hat y^0$ to $y^0$ and then
positively radiated in the $\hat z^0$ direction.
\end{example}

%ex4
\begin{example}[(The transformation model)]\label{ex4}
The transformation model (see, e.g., \citep{fra1}) provides a
paradigm for exact ancillary conditioning. A typical continuous
transformation model for a variable $y=\theta z$ has parameter $\theta$
in a smooth transformation group $G$ that
operates on an $n$-dimensional sample space for $y$;
for illustration, we assume here that the group acts coordinate by coordinate.
The natural quantile function for the $i$th coordinate is $y_i=\theta
z_i$, where $z_i$ is a coordinate reference variable with a fixed
distribution; the linear regression model with known and unknown error
scaling are simple examples.
With observed data point $y^0$, let $\hat\theta^0$ be the maximum likelihood
value and $\hat z^0$ the corresponding reference value satisfying
$y^0=\hat\theta^0 \hat z^0$. The second order approximate ancillary is
then given as $\{\theta\hat z^0\}$, which is just the usual transformation
model orbit $G \hat z^0$. If the group does not apply separately to independent
coordinates, then the present quantile approach may not be immediately
applicable; this raises issues for the construction of the trajectories
and also
for the construction of default priors (see, e.g.,
\citep{ber}). Some discussion of this in connection with
curved parameters will be reported separately. A modification achieved
by adding structure to the transformation model is given by the
structural model \citep{fra1}. This takes the reference distribution
for $z$ as the primary probability space for the model
and examines what events on that space are identifiable from an
observed response; we do not address here this alternative modelling approach.
\end{example}

%ex5
\begin{example}[(The inverted Cauchy)]
Consider a location-scale model centered at $\mu$ and scaled by $\sigma
$ with error given by the standard Cauchy; this gives the
statistical model
\[
f(y;\mu,\sigma)=\frac{1}{\curpi\sigma\{1+ (y-\mu)^2/\sigma^2\}}
\]
on the real line.
For the sampling version, this location-scale model is an example of
the transformation model discussed in the
preceding Example~\ref{ex4} and the long-accepted ancillary contour is the
half-plane (\ref{lsanc}).

McCullagh \citep{mcc2} uses linear fractional transformation results
that show that
the inversion $\tilde y=1/y$ takes the Cauchy ($\mu, \sigma$) model for $y$
into a Cauchy ($\tilde\mu, \tilde\sigma$) model for $\tilde y$, where
$\tilde\mu= {\mu}/(\mu^2+\sigma^2), \tilde\sigma= {\sigma}/(\mu
^2+\sigma^2)$.
He then notes that the usual location-scale ancillary for the derived
model does not map back to give the usual location-scale
ancillary on the initial space and would thus
typically give different inference results for the parameters; he indicates
``not that conditioning is a bad idea, but that the usual mathematical
formulation is
in some respects ad hoc and not completely satisfactory.''

We illustrate this for $n=2$ in Figure~\ref{an}.
For a data point in the upper-left portion of the plane in part (b)
for the inverted Cauchy, the observed ancillary contour is shown as a
shaded area; it is a half-plane subtended by $\mathcal{L} (1)$. When this
contour is mapped back to the initial plane in part (a), the contour
becomes three disconnected segments
with lightly shaded edges indicating the boundaries; in particular, the
line with
marks 1, 2, 3, 4, 5, 6 becomes three distinct curves again with
corresponding marks 1, 2, 3, 4, 5, 6, but
two points $(0,1), (1,0)$ on the line have no back images. Indeed, the
same type of singularity,
where a point with a zero coordinate cannot be mapped back, happens for any
sample size $n$.
Thus the proposed sample space is not one-to-one continuously
equivalent to the given sample space: points are left out and points
are created. And the quantile function used
on the proposed sample space for constructing the ancillary
does not exist on the given sample space: indeed, it is not defined at
points and is thus not continuous.

%
%f2 ###
\begin{figure}

\includegraphics{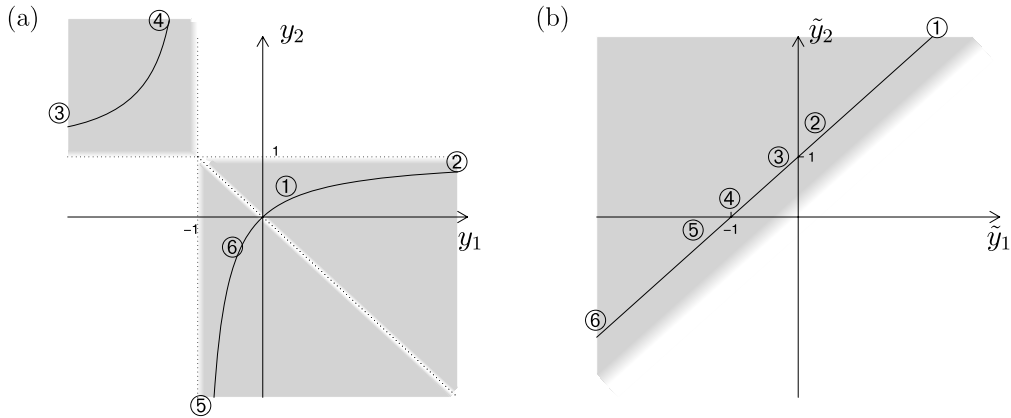}

\caption{(a) The location-scale Cauchy model for the inverted
$\tilde{y}_1 = 1/y_1$, $\tilde{y}_2 = 1/y_2$ has an ancillary contour
given by the shaded area in (b). When interpreted back for the
original $(y_1, y_2)$ the connected ancillary contour becomes three
unconnected regions, shown in (a). A line $\tilde{y}_2 = \tilde{y}_1
+1$ on the contour in (b) is mapped back to
three curved segments in (a)
and numbered points in sequence on the line are mapped back to the numbered
points on the unconnected ancillary contour.}
\label{an}
\end{figure}

The Cauchy inversion about $0$ could equally be about an arbitrary
point, say $a$,
on the real line and would lead to a corresponding ancillary. We would
thus have a wealth of competing ancillaries and a corresponding wealth
of inference procedures, and all
would have the same lack of one-to-one continuous equivalence to the
initial sample space.
While Fisher seems not
to have explicitly specified continuity as a needed ingredient for
typical ancillarity,
it also seems unlikely that he would have envisaged ancillarity without
continuity.
If continuity is included in the prescription for developing the
ancillary, then the proposed ancillary for the inverted Cauchy would
not arise.

Bayesian statistics involves full conditioning on the observed data and
familiar frequentist inference avoids, perhaps even evades,
conditioning. Ancillarity, however, represents
an intermediate or partial conditioning and, as such, offers a partial
bridging of the two extreme approaches to inference.
\end{example}

%s4 ###
\section{An asymptotic statistic}\label{sec4}
For the Normal-on-the-circle example, the exact ancillary contour was
given as the observed
contour of the radial distance $r(y_1, y_2)$: the contour is described
implicitly. By contrast, the approximate ancillary was given
as the trajectory of a point $y(\hat x^0; t)$ under change of an index
or mathematical
parameter $t$: the contour is described explicitly. For the
general context, the first approach has serious difficulties, as found
even with nonlinear
regression, and these difficulties arise with an approximate statistic
taking an approximate
value; see Example~\ref{ex2}. Accordingly, we now turn to the second, the explicit
approach, and develop the needed notation and expansions.

Consider a smooth one-dimensional contour through
some point $y_{0}$.
To describe such a contour in the implicit manner requires
$n - 1$ complementary statistics. By contrast, for the explicit method,
we write $y=y(t)$, which maps a
scalar $t$ into the sample space $\mathbb{R}^n$. More generally, for a $p$-dimensional
contour, we have $y=y(t)$ in $\mathbb{R}^n$,
where $t$
has dimension $p$ and the mapping is again into $\mathbb{R}^n$.

For such a contour, we define the row array $
V(t)= (\mathrm{d}/{\mathrm{d}t'}) y(t)=\{v_1(t),\ldots,v_p(t)\}$
of tangent vectors, where the vector
$v_\alpha(t) =(\mathrm{d}/\mathrm{d}t_\alpha)y(t)$
gives the direction or gradient of $y(t)$
with respect to change in
a coordinate $t_\alpha$.
We are interested in such a contour near a particular point $y_0 =
y(t_0)$; for convenience, we often choose $y_0$ to be the observed data
point $y^0$
and the $t_0$ to be centered
so that $t_0=0$.
In particular,
the array $V=V(t_0)$ of tangent vectors at a particular data point
$y_0$ will
be of special interest. The vectors in $V$ generate a
tangent plane $\mathcal{L} (V)$ at the point $y_0$
and this plane provides a linear
approximation to the contour. Differential geometry gives
length
properties of such vectors as the first fundamental form:
\[
V'V =
\pmatrix{
v_1\cdot v_1&\cdots& v_1\cdot v_p \cr
\vdots& & \vdots\cr
v_p\cdot v_1 & \cdots& v_p\cdot v_p
}
=
\pmatrix{
v_1^\prime v_1 & \cdots& v_1^\prime v_p\cr
\vdots& & \vdots\cr
v_p^\prime v_1 & \cdots& v_p^\prime v_p
}
;
\]
this records the matrix of inner products for the vectors $V$
as inherited from the inner product on $\mathbb{R}^n$. A
change in the parameterization $\tilde t=t(t)$ of the contour will
give different tangent vectors $V$, the same tangent plane $\mathcal{L}(V)$
and a different, but corresponding, first fundamental form.

Now, consider the derivatives of the tangents $V(t)$ at $t_0$:
\[
\label{second}
W =\frac{\mathrm{d}}{\mathrm{d}t'}V(t) \bigg|_{t=t_0}=
\pmatrix{
w_{11} &\cdots& w_{1p}\cr
\vdots & & \vdots\cr
w_{p1} & \cdots& w_{pp}
}
,
\]
where
$w_{\alpha\alpha'}=(\partial^2/\partial t_\alpha\,\partial t_{\alpha
'})y(t)|_{t=t_0}$
is an acceleration or curvature vector relative
to coordinates $t_\alpha$ and $t_{\alpha^\prime}$ at $t_0$.
We regard the array $W$ as a $p\times p$ array
of vectors in $\mathbb{R}^n$. We
could have used tensor notation,
but the approach here has the advantage
that we can write
the second degree Taylor expansion of $y(t)$ at $t_0=0$ as
%
%e4.1 ###
\begin{equation}
\label{tayloranc}
y(t)= y_0+Vt + t'Wt/2+ \cdots,
\end{equation}
which uses matrix multiplication for linearly combining the vectors in the
arrays $V$ and $W$.
Some important characteristics of the quadratic term in (\ref{tayloranc}) are
obtained by orthogonalizing the elements of $W$ to the tangent plane
$\mathcal{L} (V)$, to give residuals
\[
\tilde w_{\alpha\alpha'}= \{I-V(V'V)^{-1}V'\}w_{\alpha\alpha'} =
w_{\alpha\alpha'} -
Pw_{\alpha\alpha'};
\]
this uses the regression analysis projection matrix $P = V(V^\prime
V)^{-1}V^\prime$. The
full array $\tilde W$ of such vectors $\tilde w_{\alpha\alpha'}$
is then written
$
\tilde W = W - PW = W - VH,
$
where $H = (h_{\alpha\alpha'})$ is a $p \times p$ array of elements
$h_{\alpha\alpha'} =
(V^\prime V)^{-1}V^\prime w_{\alpha\alpha'}$; an element $h_{\alpha
\alpha'}$ is a $p
\times1$ vector, which records the regression coefficients of
$w_{\alpha\alpha'}$ on the vectors
$V$.

The array $\tilde W$ of such orthogonalized curvature vectors $\tilde
w$ is the second
fundamental form for the contour at the expansion point. Consider the
Taylor expansion
(\ref{tayloranc}) and substitute $W = \tilde W + VH$:
\begin{eqnarray*}
y(t) & = & y_0 + Vt + t^\prime(\tilde W + VH)t/2 + \cdots\\
& = & y_0 + V(t + t^\prime Ht/2) + t^\prime\tilde W t/2 + \cdots
,
\end{eqnarray*}
where we note that $t$ and $t^\prime$ are being applied to the $p
\times p$ arrays $H$
and $\tilde W$ by matrix multiplication, but the elements are $p \times
1$ vectors for $H$
and $n \times1$ vectors for $\tilde W$, and these are being combined
linearly. We can
then write $ y(t) = y_0 + V\tilde t + \tilde t \tilde W \tilde t^\prime
/2 + \cdots$
and thus have the contour expressed in terms of orthogonal curvature
vectors $\tilde w$
with the reparameterization $\tilde t = t + t^\prime H t/2 + \cdots$.
When we use this in the asymptotic setting, we will have standardized
coordinates and the
reparameterization will take the form $\tilde t = t + t^\prime H
t/2n^{1/2} + \cdots$.

%s5 ###
\section{Verifying second order ancillarity}\label{sec5}
We have used the Normal-on-the-circle example to illustrate the
proposed second order
ancillary contour $\{y(\hat x^0; t)\}$. Now, generally, let $f(y; \theta
)$ be a
statistical model with regularity and asymptotic properties as the data
dimension $ n$ increases: we
assume that the vector quantile $y(x; \theta)$ has independent scalar
coordinates and is
smooth in both the reference variable $x$ and the parameter $\theta$;
more general
conditions will be considered subsequently. For the verification, we
use a Taylor
expansion of the quantile function in terms of both $x$ and $\theta$,
and work from theory
developed in \citep{cak} and \citep{and}. The first steps involve the
re-expression of individual coordinates of $y$, $x$, and $\theta$, and
show that the
proposed contours establish a partition on the sample space; the
subsequent steps
establish the ancillarity of the contours.
\begin{longlist}[(1a)]
\item[(1a)] \textit{Standardizing the coordinates.}
Consider the statistical model in moderate deviations about $(y^0,
\hat\theta^0)$ to order $\mathrm{O}(n^{-1})$. For this, we work with coordinate
departures in
units scaled by $n^{-1/2}$. Thus, for the $i$th coordinate, we write
$y_i = \hat y_i^0 +
\tilde y_i/n^{1/2}$, $x_i = \hat x_i^0 + \tilde x_i/n^{1/2}$ and $\theta
_\alpha=
\hat\theta^0_\alpha+ \tilde\theta_\alpha/n^{1/2}$; and for
a modified $i$th quantile coordinate $\tilde y_i = \tilde y_i(\tilde
x_i, \hat\theta)$,
we Taylor expand to the second order, omit the subscripts and tildes
for temporary
clarity, and obtain $
y = x + V\theta+ (ax^{2} + 2xB\theta+ \theta'W\theta)/2n^{1/2}$,
where $V$ is the $1 \times p$ gradient of $y$ with respect to $\theta$,
$B$ is the $1\times p$ cross Hessian with respect to $x$ and $\theta$,
$W$ is the $p \times p$ Hessian with respect to $\theta$
and vector--matrix multiplication is used for combining $\theta$ with
the arrays.

\item[(1b)] \textit{Re-expressing coordinates for a nicer expansion.}
We next re-express an $x$ coordinate, writing $\tilde x = x
+ ax^2/2n^{1/2}$,
and then again omit the tildes to obtain the simpler expansion
%
%e5.1 ###
\begin{equation}
\label{scalar}
y = x + V\theta+ (2xB\theta+ \theta'W\theta)/2n^{1/2} + \cdots,
\end{equation}
to order $\mathrm{O}(n^{-1})$ for the modified $y$, $x$ and $\theta$, now in
bounded regions about
0.

\item[(1c)] \textit{Full response vector expansion.}
For the vector response $y = (y_1,\ldots, y_n)$ in quantile form, we can
compound the
preceding coordinate expansions and write
$y = x + V\theta+ (2x : B\theta+ \theta'W\theta)/ 2n^{1/2} +
\cdots,$
where $y$ and $x$ are now vectors in $\mathbb{R}^{n}$, $V = (v_1,\ldots, v_p) =
(v_\alpha)$ and $B =
(b_1,\ldots, b_p) = (b_\alpha)$ are $1 \times p$ arrays of vectors in $\mathbb{R}^{n}$,
$W = (w_{\alpha\alpha^\prime})$ is a $p\times p$ array of vectors in $\mathbb{R}^{n}$
and $x \mbox{:} B$ is a $1 \times p $ array of vectors $x \mbox{:} b$, where
the $i$th element of
the vector $x \mbox{:} b$ is the product $x_ib_i$ of the $i$th elements of
the vectors $x$
and $b$.

\item[(1d)] \textit{Eliminate the cross Hessian: scalar parameter case.}
The form of a Taylor series depends heavily on how the function and the
component
variables are expressed. For a particular coordinate of (\ref{scalar})
in (1b),
if we re-express the coordinate $ y=\tilde y + c\tilde y^2/2n^{1/2}$ in
terms of a
modified $\tilde y$, substitute it in (\ref{scalar}) and then, for
notational ease, omit the tildes, we obtain
$y+c(x+v\theta)^2/2n^{1/2}=x+v\theta+(2xb\theta+\theta^2w)/2n^{1/2}.$
To simplify this, we take the $x^2$ term over to the right-hand side
and combine it with $x$ to give a re-expressed $x$, take the $\theta x$
term over to the right-hand
side and choose $c$ so that $cv=b$ and, finally, combine the $\theta^2$
terms giving a new $w$. We then obtain
${y(x; \theta) =x+v\theta+\theta^2w/2n^{1/2}}$
with the cross Hessian removed; for this, if $v = 0$, we ignore the
coordinate as being ineffective for $\theta$.
For the full response accordingly, we then have %\label{elim2}
$
y(x; \theta) = x + v\theta+ w\theta^2/2n^{1/2}+ \cdots$
to the second order in terms of re-expressed coordinates $x$ and $y$. The
trajectory of a point $x$ is
$A(x) = \{y(x; t)\} = \{x + vt + wt^2/2n^{1/2}+\cdots\}$
to the second order as $t$ varies.

\item[(1e)] \textit{Scalar case: trajectories form a partition.}
In the standardized coordinates, the initial data point is $y^0 = 0$
with corresponding
maximum likelihood value $\hat\theta^0 = 0$; the corresponding
trajectory is
$A(0) = \{vt + wt^2/2n^{1/2} + \cdots\}.$
For a general reference value $x$, but with $\hat\theta(x) = 0$, the
trajectory is
$A(x) = \{x + vt + wt^2/2n^{1/2} +\cdots\} \nonumber= x + A(0)$.
The sets $\{A(x)\}$ with $\hat\theta(x) = 0$ are all translates of
$A(0)$ and thus
form a partition.

Consider an initial point $x_0$ with maximum likelihood value $\hat
\theta(x_0) = 0$ and let $y_1= x_0 + vt_1 + wt_1^2/2n^{1/2} + \cdots$
be a point in
the set $A(x_0) = x_0 + A(0)$.
We calculate the trajectory $A(y_1)$ of $y_1$ and show that it lies on
$A(x_0)$; the
partition property then follows and the related Jacobian effect is constant.
From the quantile function $y = x + v\theta+ w\theta^2/2n^{1/2}$, we
see that the $y$
distribution is a $\theta$-based translation of the reference
distribution described by
$x$. Thus the likelihood at $y_1$ is $l(y_1 - v\theta- w\theta
^2/2n^{1/2})$, in terms of
the log density $l(x)$ near $x_0$. It follows that $y_1 = x_0 + vt_1 +
wt_1^2/2n^{1/2}$
has maximum likelihood value $\hat\theta(y_1) = t_1$.

Now, for the trajectory about $y_1$, we calculate derivatives
\[
\frac{\mathrm{d}y}{\mathrm{d}\theta} =v+w\theta/n^{1/2}, \qquad \frac{\mathrm{d}^2y}{\mathrm{d}\theta^2}=w/n^{1/2},
\]
which, at the point $y_1=vt_1 + wt^2_1/2n^{1/2}$ with $\theta=\hat\theta
(y_1)$, gives
\[
V(y_1) =v+wt_1/n^{1/2}, \qquad  W(y_1)=w/n^{1/2},
\]
to order $\mathrm{O}(n^{-1})$.
We thus
obtain the trajectory of the point $y_1$:
\begin{eqnarray*}
A(y_1) &=& \{x_0 + vt_1+
wt_1^2/2n^{1/2}+(v+wt_1/n^{1/2})t+wt^2/2n^{1/2}\} \\
&=& \{x_0 + vT +wT^2/2n^{1/2}\}
\end{eqnarray*}
under variation in $t$.
However, with $T=t_1+ t$, we have just an arbitrary point on the
initial trajectory. Thus the mapping $y \to A(y)$ is well defined and
the trajectories generate a partition, to second order in moderate
derivations in
$\mathbb{R}^n$. In the standardized coordinates, the Jacobian effect is constant.

\item[(1f)] \textit{Vector case: trajectories form a partition}.
For the vector parameter case, we again use
standardized coordinates and
choose a parameterization that
gives orthogonal curvature vectors $w$
at the observed data point $y^0$. We then
examine scalar parameter change on some line through $\hat\theta(y^0)$.
For this, the results above give a trajectory and
any point on it reproduces the trajectory under that scalar parameter.
Orthogonality ensures that the vector maximum likelihood value is on
the same line
just considered. These trajectories are,
of course, part of the surface defined by $\{Vt+t'Wt/2n^{1/2}\}$.
We then use the partition property of the individual trajectories
as these apply perpendicular to the surface; the surfaces are thus part
of a partition. We can then write the trajectory of a point $x$ as a set
%
%e5.2 ###
\begin{equation}
\label{taylor3}
A(x)=\{ x + Vt + t'Wt/ 2n^{1/2} + \cdots\dvt t\}=x + A(0)
\end{equation}
in a partition to the second order in moderate deviations.

\item[(2a)] \textit{Observed information standardization}.
With moderate regularity, and
following \citep{fra5} and \citep{fra10}, we have a limiting Normal
distribution conditionally on
$ y^0+\mathcal{L}(V)$. We then rescale the parameter
at $\hat\theta^0$ to give identity observed information and thus an identity
variance matrix for the Normal distribution to second order.
We also have a limiting
Normal distribution conditionally on $ y^0+\mathcal{L}(V,W)$; for this,
we linearly modify the vectors in $W$ by rescaling and regressing on
$\mathcal{L}(V)$
to give distributional orthogonality to $\hat\theta$ and identity
conditional variance
matrix to second order.

\item[(2b)] \textit{The trajectories are ancillary: first derivative parameter change}.
We saw in the preceding section that key local properties of a statistical
model were summarized by the tangent vectors $V$ and the curvature vectors
$W$, and that the latter can, to advantage, be taken to be orthogonal
to the
tangent vectors. These vectors give local coordinates for the model and
can be replaced by
an appropriate subset if linear dependencies are present.

First, consider the conditional model given the
directions corresponding to the span $y^0+\mathcal{L}\{V, W\}$. From the ancillary
expansion (\ref{taylor3}), we have that change of $\theta$ to the
second order
moves points within the linear space $y^0+\mathcal{L}\{V, W\}$; accordingly,
this conditioning is ancillary.
Then, consider the further conditioning to an alleged ancillary
contour, as described by (\ref{taylor3}). Also, let $y_0$
be a typical point having
$\hat\theta(y_0)=\hat\theta^0$ as the corresponding maximum likelihood value;
$y_0$ is thus on the observed maximum likelihood contour.

Now, consider a rotationally symmetric Normal distribution on the $(x,
y)$ plane
with mean $\theta$ on the $x$ axis and let $a=y+cx^2/2$ be linear in $y$
with a quadratic adjustment with respect to~$x$. Then $a=a(x,y)$ is
first-derivative
ancillary at $\theta=0$.
For this, we assume,
without loss of generality, that the standard deviations
are unity. The marginal density for $a$ is then
\[
f(a;\theta)=\int_{-\infty}^{\infty}\phi(x-\theta)\phi(a-cx^2/2)\,\mathrm{d}x ,
\]
which is symmetric in $\theta$; thus
$(\mathrm{d} / \mathrm{d}\theta)f(a;\theta) |_{\theta=0}=0,$
showing
that the distribution of $a$ is first-derivative ancillary
at $\theta=0$ or, more intuitively, that the amount of probability on
a contour of $a$ is first-derivative free of $\theta$ at $\theta=0$.
Of course, for this, the $y$-spacing between contours of $a$ is constant.

Now, more generally, consider an asymptotic distribution for $(x,y)$
that is first order
rotationally symmetric Normal with mean $\theta$ on the $y=0$ plane;
this allows
$\mathrm{O}(n^{-1/2})$ cubic contributions. Also, consider an $s$-dimensional variable
$a=y+Q(x)/2n^{1/2}$ which is a quadratic adjustment of $y$. The
preceding argument
extends to show that $a(y)$ is first-derivative ancillary: the two $\mathrm{O}(n^{-1/2})$
effects are zero and the combination is of the next order.

\item[(2c)] \textit{Trajectories are ancillary: parameter change in moderate deviations.}
Now, consider a statistical model $f(y;\theta)$ with data point $y^0$ and
assume regularity, asymptotics and smoothness of the quantile
functions. We examine the parameter trajectory $\{y(\hat x^0; t)\}$
in moderate deviations under change in $t$. From the preceding
paragraph, we then have
first-derivative
ancillarity at $\theta=\hat\theta=0$. But this holds for each expansion
in moderate
deviations
and we thus have ancillarity in moderate deviations. The key here
has been to use the expansion form about the point that has $\hat\theta$
equal to the parameter value being examined.
\end{longlist}
%s6 ###
\section{Discussion}\label{sec6}

\begin{longlist}[(iii)]
\item[(i)] \textit{On ancillarity.}
The \hyperref[intro]{Introduction} gave a brief background on ancillary statistics
and noted that an ancillary is typically viewed as a statistic with a
parameter-free distribution; for some recent discussion, see \citep{fra4}.
Much of the literature is concerned with difficulties that can arise
using this third
Fisher concept, third after sufficiency and likelihood: that maximizing
power given size
typically means not conditioning on an ancillary; that shorter
on-average confidence intervals typically mean ignoring ancillary
conditioning; that techniques that are conditional on an ancillary are
often inadmissible; and more. Some of the difficulty may hinge on whether
there is merit in the various optimality criteria themselves. However,
little in the literature seems
focused on the continued evolution and development of this Fisher
concept, that is, on what
modifications or evolution can continue the exploration initiated in
Fisher's original papers (see \citep{fis1,fis2,fis3}).

\item[(ii)] \textit{On simulations for the conditional model.}
The second order ancillary in moderate deviations has contours that
form a partition, as shown in the preceding section. In the modified or
re-expressed coordinates, the contours
are in a location relationship and, correspondingly, the Jacobian
effect needed for the
conditional distribution is constant. However, in the original coordinates,
the Jacobian effect would typically not be constant and its effect
would be needed
for simulations. If the parameter is scalar, then the effect is
available to the second order
through the divergence function of a vector field; for some discussion
and examples, see \citep{fra2}. For a vector parameter, generalizations
can be implemented, but we do not
pursue these here.

\item[(iii)] \textit{Marginal or conditional.}
When sampling from a scalar distribution
having variable $y$ and moderate regularity,
the familiar central limit theorem gives a limiting Normal distribution
for the sample average $\bar y$ or
sample sum $\sum y_i$.
From a geometric view, we have probability in
$n$-space and contours determined by
$\bar y$, contours that are planes perpendicular to the $1$-vector.
If we then collect the probability on a contour, plus or
minus a differential, and deposit it, say, on the intersection
of the contour with the span
$\mathcal{L} (1)$ of the $1$-vector, then we obtain a limiting Normal
distribution on $\mathcal{L} (1)$, using
$\bar y$ or $\sum y_i$ for location on that line.

A far less familiar Normal limit result applies in the same general
context, but with a totally different geometric decomposition. Consider lines
parallel to the $1$-vector, the affine cosets of $\mathcal{L} (1)$. On these
lines, plus or minus a differential, we then obtain a limiting Normal
distribution for location say $\bar y$ or $\sum y_i$.
In many ways, this conditional, rather than marginal, analysis
is much stronger and more useful.
The geometry, however, is different,
with planes perpendicular to $\mathcal{L} (1)$ being replaced
by points on lines parallel to $\mathcal{L} (1)$.

This generalizes giving a limiting conditional Normal distribution
on almost arbitrary smooth contours in a partition and it has wide application
in recent likelihood inference theory. It also provides third order accuracy
rather than the first order accuracy associated with the usual geometry.
In a simple sense, planes are replaced by lines or
by generalized contours and much stronger, though
less familiar, results are obtained.
For some background based on Taylor expansions of log-statistical
models, see \citep{cak,che} and \citep{and}.
\end{longlist}

\section*{Acknowledgements}

This research was supported by the Natural Sciences and Engineering
Research Council of
Canada. The authors wish to express deep appreciation to the referee
for very incisive comments. We also offer special thanks to Kexin Ji
for many contributions and support with the manuscript and the diagrams.

\printhistory

\end{document}